\documentclass{article}
\usepackage{graphics}
\usepackage{amsfonts}

\newcommand{\R}{\mathbb{R}}

\newcommand{\fg}{\mathfrak{g}}

\newcommand{\cA}{\mathcal{A}}
\newcommand{\cC}{\mathcal{C}}
\newcommand{\cG}{\mathcal{G}}

\newcommand{\hoopg}{\mathcal{HG}}
\newcommand{\loopg}{\mathcal{LG}}
\newcommand{\hoopgw}{\mathcal{HG}^\omega}
\newcommand{\loopgw}{\mathcal{LG}^\omega}

\newcommand{\clos}{\overline}

\newcommand{\Ker}{\mathop{\mathrm{Ker}}\nolimits}
\newcommand{\Aut}{\mathop{\mathrm{Aut}}\nolimits}
\newcommand{\range}{\mathop{\mathrm{range}}\nolimits}

\newcommand{\mathrlap}[2]{%
    \hbox to \mathrlapwidth {$#1{{}#2}\mathsurround=0pt$\hss}}
\newcommand{\mrlap}[1]{\def\mathrlapwidth{#1}\mathpalette\mathrlap}

\newenvironment{eqnalign*}%
{\displaymath\def\eqcr{\cr\noalign{\vskip 4pt}\relax} \let\\\eqcr
\vbox\bgroup \halign\bgroup \hfill
$\displaystyle{##{}}$&$\displaystyle{{}##}$\hfill\cr\relax}%
{\crcr\egroup\unskip\unskip\unskip\egroup \enddisplaymath}

\newtheorem{teo}{Theorem}[section]
\newtheorem{lem}[teo]{Lemma}
\newtheorem{prop}[teo]{Proposition}
\newtheorem{cor}[teo]{Corollary}

\newtheorem{defi}[teo]{Definition}

\newcommand{\Proof}{\par\vskip 5pt plus 2pt minus 1pt\noindent {\it Proof: }}
\newcommand{\ProofP}[1]{\par\vskip 5pt plus 2pt minus 1pt\noindent
                        {\it Proof #1:}}
\newcommand{\QED}{\unskip\unskip\unskip\penalty 0\hbox to 0pt {}
                  \nobreak\hfil QED{\parfillskip=0pt\par}%
                  \vskip 6pt plus 2pt minus 1pt\relax}

\begin{document}

\title{Groups of loops and hoops}

\author{Pablo Spallanzani\\
	{\tt pablo@cmat.edu.uy}\\
        Centro de Matem\'atica, Facultad de Ciencias\\
        Igua 4225, Montevideo CP11400, Uruguay}

\date{11/8/1999}

\maketitle

\begin{abstract}
The approaches to quantum field theories based in the so called
loop representation deserved much attention recently.
In it, closed curves and holonomies around them play a central role.
In this framework the group of loops and the group of hoops have been defined,
the first one consisting in closed curves quotient with the equivalence 
relation that identifies curves differing in retraced segments, and the second
one consisting in closed curves quotient with the equivalence relation that
identifies curves having the same holonomy for every connection in a fiber
bundle.
The purpose of this paper is to clarify the relation between hoops and loops,
or in other words, to give a description of the class of holonomy equivalent
curves.
\end{abstract}

\section{Introduction}

An important step in the construction of quantum field theories is the
definition of the space of states $L^2(\cA/\cG)$.
This is done in \cite{AL1,BS1} by first constructing generalized measures on
$\cA/\cG$.
In these constructions the notions of group of loops, holonomy around
loop and group of hoops play a central role (the precise definitions are
stated below).

Given a differentiable manifold $M$ and a point $o$ of $M$ we construct the
space of closed curves in $M$, $\Omega$, as the set of piecewise regular curves
$\alpha:[0,1] \to M$. A curve $\beta:[a,b] \to M$ is regular if there exists
$\epsilon>0$ and a differentiable (or analytic) curve
$\gamma:(a-\epsilon, b+\epsilon) \to M$ such that $\beta$ and $\gamma$
coincide in $[a,b]$, and we say that a curve $\alpha:[0,1] \to M$ is
piecewise regular if exists a partition of $[0,1]$, $0=t_0<t_1<\cdots<t_n=1$
such that $\alpha$ restricted to each of the intervals $[t_{i-1},t_i]$ is
regular.

In $\Omega$ we can define the inverse of a curve $\alpha^{-1}(t)=\alpha(1-t)$,
and the composition of curves
\[
    \alpha\beta(t)=\left\{\begin{array}{ll}
                            \alpha(2t)  & \hbox{if } t <   1/2 \\
                            \beta(2t-1) & \hbox{if } t \ge 1/2
                          \end{array} \right.
\]
if we identify curves that only differ in a reparameterization the composition
is an associative operation, but in general $\alpha\alpha^{-1} \ne c$ 
($c$ being the constant curve), to make $\Omega$ in a group we need to 
introduce a further equivalence relation.
One possibility is to identify curves differing in retraced segments, that is
we identify $\alpha\beta$ with $\alpha\rho\rho^{-1}\beta$, the group
obtained is called the group of loops and is denoted by $\loopg$ or
$\loopgw$ if we work with analytic curves.
Other possible identification is, given a principal bundle $(E,M,G,\pi)$, $G$ a
Lie group, identify two curves $\alpha$ and $\beta$ if they have the same
holonomy for every connection in the bundle, the group obtained this way is
called the group of hoops and is denoted by $\hoopg$ or $\hoopgw$ in the
analytic case.
The purpose of this paper is to clarify the relation between $\loopg$ and
$\hoopg$ and how $\hoopg$ depends on the Lie group $G$.
In particular we obtain results for piecewise differentiable loops without
making any assumptions on the Lie group.

Consider the infinite set of symbols $e_1,e_2,\ldots$ and
$e_1^{-1},e_2^{-1},\ldots$, and let $E$ be the set of words in that symbols
including the null word, a word is a finite ordered list of symbols (ex.\@
$e_3e_1^{-1}e_2$ is a word).
If we define the product of word as the concatenation and identify words that
differ by ``canceling opposite symbols'', that is
$w_1e_ie_i^{-1}w_2 \sim w_1w_2$, then $E$ is a free group.
Let us define $E_G$, the group of {\em identities} of $G$, as the subgroup of
$E$ consisting in words, say $e_2e_3e_1^{-1}$, such that if we assign to each
$e_i$ a element $g_i$ of $G$ and multiply these in the way specified by the
word (as $g_2g_3g_1^{-1}$) the result is the identity of $G$ no matter what
choice of $g_i$ (ex.\ if $G$ is abelian $e_1e_2e_1^{-1}e_2^{-1}$ is an
identity). In other words, $E_G$ is the intersection of the kernels of every
homomorphism of groups from $E$ to $G$
\[
    E_G = \bigcap_{f \in \mrlap{2em}{\hom{(E,G)}}} \ker f.
\]
Now we define $E_G(\loopg)$ as the subgroup of $\loopg$ generated by the loops
obtained in the following way: for every word in $E_G$ (such as
$e_2e_1e_3^{-1}$) and every assignment of a loop $\alpha_i$ to each of the
symbols $e_i$ take the product of $\alpha_i$ in the same way as the word
(ex.\ $\alpha_2\alpha_1\alpha_3^{-1}$). Or equivalently
\[
    E_G(\loopg) = \bigcup_{f \in \mrlap{2em}{\hom{(E,\loopg)}}} f(E_G).
\]

Now we can state the main results, first in the analytic case.

\begin{teo}\label{teo:hg-an-loop}
For $G$ a connected Lie group, $\hoopgw = \loopgw/E_G(\loopgw)$.
\end{teo}

This result is complemented with results about $E_G$ of
section~\ref{sec:ident-groups}.

\begin{teo}
If $G$ is abelian then $E_G$ is generated by elements of the form
$e_ie_je_i^{-1}e_j^{-1}$.
\end{teo}

\begin{teo}
If $G$ is connected and non solvable then has no non trivial identities.
\end{teo}

Then we have the following corollaries (see~\cite{AL1}).

\begin{cor}
If $G$ is abelian then $\hoopgw = \loopgw/[\loopgw,\loopgw]$.
\end{cor}

\begin{cor}\label{cor:hg-an-loop-nonsolv}
If $G$ is connected and non solvable then $\hoopgw = \loopgw$.
\end{cor}

In section~\ref{sec:diff-loop} we show through a example that 
theorem~\ref{teo:hg-an-loop}
is not valid in the differentiable case. However we have:

\begin{teo}\label{teo:hg-diff-loop}
For $G$ a connected Lie group, $\hoopg = \loopg/\clos{E_G(\loopg)}$.
\end{teo}

Where $\clos{E_G(\loopg)}$ is the closure of $E_G(\loopg)$ in the
quotient topology arising from the $C^N$ topology of curves for any $N$.
The topology of the loop space is discussed in
section~\ref{sec:topology-loop} where we also show that the topology
introduced by Barret~\cite{Ba} coincides with the usual $C^N$
topology.
However we can generalize corollary~\ref{cor:hg-an-loop-nonsolv}
for a non solvable group even in the case of piecewise differentiable loops,
as the following theorem shows.

\begin{teo}\label{teo:hg-diff-loop-nonsolv}
If $G$ is connected and non solvable then $\hoopg = \loopg$.
\end{teo}

\section{Analytic loops}\label{sec:an-loop}

First we consider the case of analytic loops, this case is simpler because of
the way in which analytic curves intersect, they either intersect in finitely
many points or they intersect in a segment.
From this we obtain a decomposition of a loop in independent loops.

Let us define what we mean by {\em independent loops}, we say that a loop
$\alpha$ has a segment $\rho$ that is traced once if there exists curves
$\beta$ and $\gamma$ such that $\alpha=\beta\rho\gamma$ and $\beta$ and
$\gamma$ don't intersect $\rho$ except at the endpoints. A set of loops
$\alpha_1,\ldots,\alpha_n$ is independent if each loop $\alpha_i$ has a
segment $\rho_i$ traced once and the segments $\rho_i$ do not intersect.

\begin{teo}\label{teo:decomp-an-loop}
Every loop can be decomposed in product of independent loops.
\end{teo}

\Proof
The loop $\gamma$ is piecewise analytic thus it can be written as a product of
analytic curves $\gamma=\rho_1\ldots\rho_n$, the curves $\rho_i$
intersect each other in finitely many points or in a common segment thus each
$\rho_i$ can be decomposed in segments that intersect only at the endpoints
or coincide, then $\gamma=\alpha_{i_1}^{s_1}\ldots\alpha_{i_k}^{s_k}$, where
$s_j$ is either $1$ or $-1$.

Let us denote by $e_-(\alpha_i)$ the initial point of $\alpha_i$ and
$e_+(\alpha_i)$ the final point. Let $E$ denote the set of all the endpoints
of all $\alpha_i$, for each point $p\in E$ chose a curve $\beta(p)$ from $o$ to
$p$ that does not intersect the segments $\alpha_i$, and let
$\gamma_i=\beta(e_-(\alpha_i))\alpha_i\beta(e_+(\alpha_i))$ then
$\gamma=\gamma_{i_1}^{e_1}\ldots\gamma_{i_k}^{e_k}$.
\QED

\begin{lem}\label{lem:chose-holonomy}
Let $(E,M,G,\pi)$ be a principal bundle with $G$ a connected Lie group and
$\alpha$ a loop in $M$ with a segment traced once, chose any element $g$ of
$G$, then there is a connection $\theta$ in $E$ such that
$H_\theta(\alpha)=g$.
\end{lem}

\Proof
The loop $\alpha$ has a segment traced once, thus we can find a local
parameterization of $M$ such that
\begin{enumerate}
\item its domain contains $I=[0,1]^n$, in what follows we identify points in
      $I$ with its images in $M$ and we fix a trivialization of the bundle over
      $I$.
\item $\alpha=\beta\gamma\xi$ such that $\beta$ and $\xi$ don't have points in
      $I$ except its endpoints.
\item $\gamma$ is the segment from $a=(0,1/2,\ldots,1/2)$ to
      $b=(1,1/2,\ldots,1/2)$ in $I$.
\end{enumerate}

\begin{figure}[h]
\centering
\includegraphics{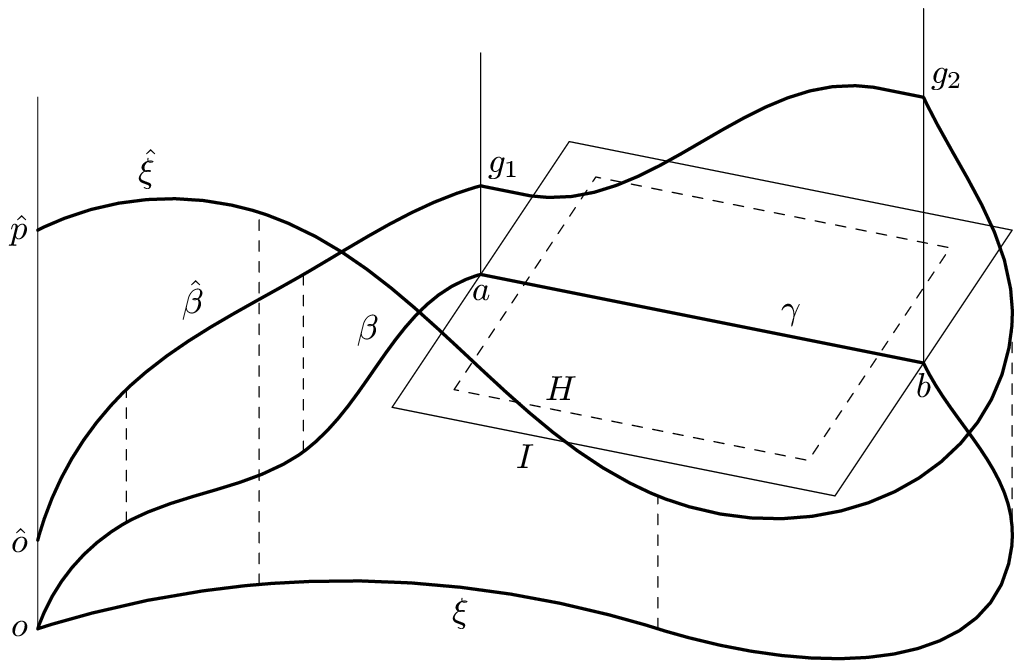}
\caption{}\label{fig:conshol}
\end{figure}

Let $A$ be any connection, take a small deformation of $A$ such that $A$ is
flat over $I$.
Take $\hat{o}$ and $\hat{p}$ in the fiber over $o$ such that
$\hat{o}=\hat{p}g$, let $\hat{\beta}$ the horizontal lift of $\beta$ that
starts in $\hat{o}$ and $\hat{\xi}$ the horizontal lift of $\xi$ that ends in
$\hat{p}$.
Next fix a trivialization of the bundle over $I$, thus elements in fibers over
points of $I$ can be identified with elements of $G$. Let $g_1$ be the
endpoint of $\hat\beta$ in the fiber over $a$ and $g_2$ the endpoint of
$\hat\xi$ in the fiber over $b$, see figure, $G$ is connected then there
exists a curve $s:[0,1] \to G$ such that $s(t)=g_1$ for $t \in [0, \epsilon)$
and $s(t)=g_2$ for $t$ in $(1-\epsilon,1]$ and let $\rho:[0,1] \to \R$
differentiable such that $\rho(t)=1$ for $t\in [\epsilon,1-\epsilon]$ and
$g(t)=0$ for $t \in [0,\epsilon/2) \cup (1-\epsilon/2,1]$,
take the connection 
$A'_x=L_{s(x_1)^{-1}}\dot{s}(x_1)\rho(x_2)\ldots\rho(x_n)dx_1$,
let $B$ the connection defined as $A$ outside $I$ and as $A'$ in $I$ 
(note that this is a smooth connection because $A$ is flat in $I$ and the way 
in which $s$ and $\rho$ were chosen), then $H_B(\alpha)=g$.
\QED

Note that the proof of this lemma requires to change a given connection only
in a small neighborhood of a point in the segment traced once, then we can
use it to prove the next proposition.

\begin{prop}\label{prop:chose-holonomy2}
Let $(E,M,G,\Pi)$ be a principal bundle with $G$ a connected Lie group and
$\alpha_1,\ldots,\alpha_n$ be independent loops then for every
$(g_1,\ldots,g_n)$ in $G^n$ there is a connection $A$ such that
$H_A(\alpha_i)=g_i$ for $i=1,\ldots,n$.
\end{prop}

Now we state and prove the main theorem of the section.

\begin{teo}
For $G$ a connected Lie group, $\hoopgw=\loopgw/E_G(\loopgw)$.
\end{teo}

\Proof
Let $\alpha$ be a loop in $E_G(\loopgw)$ then there exist a word in $E_G$, for
example $e_1e_2e_1^{-1}e_3e_2^{-1}$, such that
$\alpha=\alpha_1\alpha_2\alpha_1^{-1}\alpha_3\alpha_2^{-1}$ then
for every connection $A$, if we define $g_i=H_A(\alpha_i)$, the holonomy of
$\alpha$ is $H_A(\alpha)=g_1g_2g_1^{-1}g_3g_2^{-1}=e$.
Conversely if $\alpha$ is a loop not in $E_G(\loopgw)$ then by
theorem~\ref{teo:decomp-an-loop} there exists independent loops
$\alpha_1,\ldots,\alpha_n$ such that
$\alpha=\alpha_{i_1}^{s_1}\ldots\alpha_{i_k}^{s_k}$.
Note that $e_{i_1}^{s_1} \ldots e_{i_k}^{s_k}$ is not an identity of $G$
because $\alpha\not\in E_G(\loopgw)$, then there exist $g_1,\ldots,g_n$ in
$G$ such that $g_{i_1}^{s_1} \ldots g_{i_k}^{s_k}\ne e$; by
proposition~\ref{prop:chose-holonomy2} there exists a connection $A$ such that
$H_A(\alpha_i)=g_i$ then $H_A(\alpha)\ne e$.
\QED

\section{Identities in Lie groups}\label{sec:ident-groups}

In this section we prove the following theorems

\begin{teo}\label{teo:ident-abelian}
If $G$ is abelian then $E_G$ is generated by elements of the form
$e_ie_je_i^{-1}e_j^{-1}$.
\end{teo}

\begin{teo}\label{teo:ident-nonsolv}
If $G$ is connected and non solvable then it has no non trivial identities.
\end{teo}

\ProofP{of~\ref{teo:ident-abelian}}
If $g$ is abelian then $E_G$ contains all words of the form
$e_ie_je_i^{-1}e_j^{-1}$.
Conversely, if $e_{i_1}^{s_1}\ldots e_{i_k}^{s_k}$ is an identity all the
words formed by reordering terms are also identities because $G$ is abelian,
thus is sufficient to prove that $e_1^{a_1}\ldots e_n^{a_n}$ is an identity
iff $a_i=0$ for $i=1,\ldots,n$.
If $a_j=m \ne 0$ then take $g$ an element of $G$ such that $g^m \ne e$ (for
example if $v$ is a vector in $\fg$ such that $\exp v \ne e$ the take
$g=\exp v/m$), define $g_i=e$ if $i\ne j$ and $g_j=g$ then
$g_1^{a_1}\ldots g_n^{a_n}=g^m \ne e$ thus $e_1^{a_1}\ldots e_n^{a_n}$ is not
an identity of $G$.
\QED

To prove~\ref{teo:ident-nonsolv} we show that if $G$ is non solvable then for
every $n$ it has a free subgroup with $n$ generators. We will use the
following theorem due to Tits~\cite{Ti1}.

\begin{teo}\label{teo:teo-tits}
Let $G \subset GL(V)$ be a subgroup, $V$ a finite-dimensional vector space
over a field of characteristic $0$. Then $G$ has a free subgroup with $n$
generators for every $n$ or $G$ has a solvable subgroup of finite index.
\end{teo}

First we recall the definition of solvable groups. Let $G$ be a group.
The derived group $G'$ is the subgroup of $G$ generated by elements of the
form $xyx^{-1}y^{-1}$, $x,y\in G$, then define by induction
\[
        G^{(1)}=G'   \quad\quad  G^{(n+1)}=\big(G^{(n)}\big)',
\]
then $G$ is {\em solvable} if $G^{(n)}$ is the trivial group for some $n$.
Next we prove the following theorem.

\begin{teo}\label{teo:ident-nonsolv1}
If $G$ is a subgroup of $GL(n,\R)$ connected and non solvable then for every
$n$ it has a free subgroup with $n$ generators.
\end{teo}

\Proof
Suppose that $G$ does not contain a free subgroup with $n$ generators, then by
theorem~\ref{teo:teo-tits} $G$ has a solvable subgroup of finite index $H$.
Let $\clos{H}$ be the closure of $H$ in $G$; then $\clos{H}$ is a solvable
subgroup of finite index of $G$, then either $\clos{H}=G$ which is absurd
because $G$ is nonsolvable, or the index of $\clos{H}$ in $G$ is greater than
$1$ then $G$ is union of finitely many closed disjoint subset (the cosets
of $\clos{H}$) which is absurd because $G$ is connected.
\QED

To prove theorem~\ref{teo:ident-nonsolv} we use the adjoint representation of
a Lie group, $Ad:G \to \Aut(\fg)$, $Ad(g)v=da_g v$ where $a_g :G \to G$,
$a_g(x)=gxg^{-1}$.

\ProofP{of~\ref{teo:ident-nonsolv}}
We need to show that $Ad(G)$ is a connected nonsolvable subgroup of
$\Aut(\fg)$. Clearly $Ad(G)$ is connected because $Ad$ is continuous.
Suppose that $Ad(G)$ is solvable, that is, there exist $n$ such that
$Ad(G)^{(n)}=\{e\}$, but $Ad(G)^{(n)}=Ad(G^{(n)})$ then
$G^{(n)} \subset \Ker Ad=Z(G)$ ($Z(G)$ is the set of all elements in $G$
that commute with every other element of $G$),
then $G^{(n+1)}=\{e\}$ which is absurd because $G$ is nonsolvable.
Then $Ad(G)$ is a connected nonsolvable subgroup of $\Aut(\fg)$ and by
theorem~\ref{teo:ident-nonsolv1}, $Ad(G)$ has no non trivial identities thus
$G$ has no non trivial identities.
\QED

\section{Topology of the loop space}\label{sec:topology-loop}

In this section we discuss several ways to give a topology to the loop space,
we work in the space of parameterized paths in $M$, let
\[
    P^N=\{ \gamma:[0,1] \to M: \hbox{$\gamma$ is piecewise $C^N$} \}
\]
and we define $P^\infty=\bigcap_{N>0} P^N$.
We define the $C^N$ topology in $P^N$ giving a subbase of open sets.
Let $\phi:U\subset M \to \R^d$ be a coordinate system in $M$, $a<b \in [0,1]$
and $\gamma$ a curve such that $\gamma|_{[a,b]} \subset U$, then we define
\begin{eqnalign*}
    U^N_{\phi,a,b}(\gamma,\epsilon) =
        \{ \alpha \in P^N : \alpha|_{[a,b]} \subset U \hbox{ and }
            |\alpha^{(n)}(x)-&\gamma^{(n)}(x)|<\epsilon,            \\
            &\forall n \le N, x \in [a,b]  \}
\end{eqnalign*}
where we identify $\alpha$ with $\phi\circ\alpha$ and we denote by
$\alpha^{(n)}$ the $n$-th derivate of $\alpha$.
We take this family of sets with $N$ fixed as a subbase of the $C^N$
topology in $P^N$ and take the family of this sets for all $N$ as a subbase
of the $C^\infty$ topology in $P^\infty$.

We now give another characterization of this topology.
Define a $C^N$ homotopy, where possibly $N=\infty$,
as function $\Phi:U \to P^N$ that is obtained from a $C^N$ function
$\phi:U\times [0,1] \to M$ with $U$ a open set of $\R^n$, 
the finest topology in which all the $C^N$ homotopies are continuous
is the Barret~\cite{Ba} topology.
We affirm that the Barret topology coincides with the $C^N$ topology.
We give a proof for the $C^\infty$ case, the $C^N$ case is
similar.

In what follows we consider that all the curves are contained in the domain
of a coordinate system (if not we can divide paths in smaller pieces) and
thus we identify them with paths in $\R^d$.
As obviously all the $C^\infty$ homotopies are continuous in the $C^\infty$ 
topology, closed sets in the $C^\infty$ topology are closed in the topology 
generated by $C^\infty$ homotopies, the converse follows from the following
lemma.

\begin{lem}
If $\gamma_n$ is a sequence of $C^\infty$ curves in $\R^n$ converging in the
$C^\infty$ topology to $\gamma$ then there is a homotopy
$\Phi:(-1,1) \to P^\infty$ such that $\Phi(0)=\gamma$ and
$\alpha_n=\Phi(2^{-n})$ is a subsequence of $\gamma_n$.
\end{lem}

\Proof
Take $\rho:[0,1] \to [0,1]$ a $C^\infty$ function such that $\rho(0)=0$,
$\rho(1)=1$, $\rho^{(n)}(0)=\rho^{(n)}(1)=0$ for all $n \ge 0$ and let
$a_n = \max_{x\in [0,1], k \le n} |\rho^{(k)}(x)|$.
Take $\alpha_n$ a subsequence of $\gamma_n$ such that
$\alpha_n \in U^N(\gamma,2^{-N^2-N-1}a_N^{-1})$ for all $n \ge N-1$.
Define $\phi:(-1,1)\times[0,1] \to \R^d$ as
\[
    \phi(s,t)=\left\{ \begin{array}{ll}
           \gamma(t) & \hbox{if } s \le 0 \\
           (1-\rho(2^n s - 1))\alpha_n(t) + \rho(2^n s - 1)\alpha_{n-1}(t)
                     & \hbox{if } 2^{-n} < s \le 2^{-n+1}
    \end{array} \right.
\]
Then obviously $\phi(s,t)$ is $C^\infty$ when $s \ne 0$, we have to show that
\[
    \frac{\partial^{k+l}\phi}{\partial s^k \partial t^l} \to 0
\]
when $s \to 0$, take $n>k+l$, then for $2^{-n} < s \le 2^{-n+1}$ we have
\begin{eqnalign*}
    \left|\frac{\partial^{k+l}\phi}{\partial s^k \partial t^l}\right| &=
        |-2^{nk}\rho^{(k)}(2^n s - 1)\alpha_n^{(l)}(t) +
        2^{nk}\rho^{(k)}(2^n s - 1)\alpha_{n-1}^{(l)}(t)| \\
    &\le 2^{n^2} a_n |\alpha_n^{(l)}(t) - \alpha_{n-1}^{(l)}(t)|
        \le 2^{n^2} a_n 2^{-n^2-n} a_n^{-1} \le 2^{-n}.
\end{eqnalign*}
Then $\phi(s,t)$ is a $C^\infty$ function thus $\Phi$ is a $C^\infty$
holonomy such that $\alpha_n=\Phi(2^{-n})$ and $\gamma=\Phi(0)$.
\QED

\section{Differentiable loops}\label{sec:diff-loop}

The key of the proof of theorem~\ref{teo:hg-an-loop} for analytic loops
was the theorem of decomposition of a loop in product of independent loops
(theorem~\ref{teo:decomp-an-loop}), in the case of differentiable loops this
theorem is not valid because differentiable curves can intersect in
complicated ways.

\begin{figure}[h]
\centering
\includegraphics{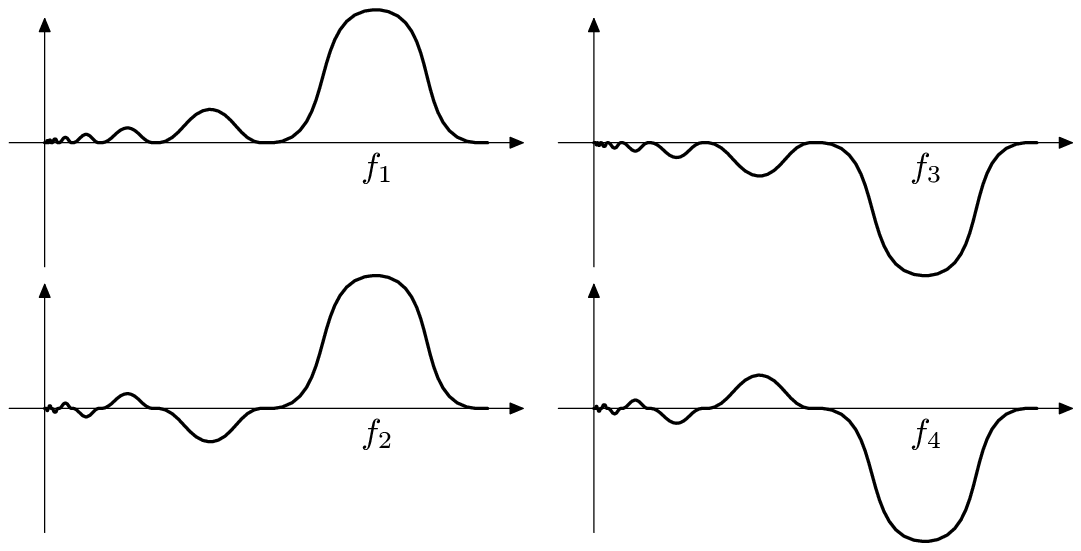}
\caption{}\label{fig:example}
\end{figure}

For example let $\rho:[0,1] \to [0,1]$ a differentiable function such that
$\rho^{(n)}(0)=\rho^{(n)}(1)=0$ for all $n \ge 0$ and let
$a_n = \max_{x\in [0,1], n>0} |\rho^{(n)}(x)|$.
Define $f_1, f_2 : [0,1] \to \R$ as
\[
    f_i(x)= (3-2i)^n \frac{1}{2^n a_n} \rho(2^n x -1) \hbox{ if }
        2^{-n}< x \le 2^{-n+1}, \quad i=1,2
\]
and $f_3=-f_1$, $f_4=-f_2$, see figure~\ref{fig:example}.

Take the curves $c_i(x)=(x,f_i(x))$ and $c=c_1 c_2^{-1} c_3 c_4^{-1}$
then $c$ has trivial holonomy for any connection in a bundle with abelian
structure group $G$ but it is not in $E_G(\loopg)$.

However we have another decomposition available~\cite{BS1}.

Let us begin with some definitions, let $T$ be a family of curves
$\alpha_1,\ldots,\alpha_n$ in $M$, $\alpha_i:[0,1] \to M$, we define
$\range(T)$ as the union of the images of the curves in $M$,
a point $p$ of $\range(T)$ is a {\em regular point} if there is a
neighborhood $U$ of $p$ such that $U \cap \range(T)$ is an embedded segment in
$M$.

\begin{defi}
Let $p$ be a regular point, the {\em type} of $p$ is the subgroup of
$G^n$ generated by the elements in which the $i$-th component is $e$ 
if $\alpha_i$ do not pass through $p$ and the other components are all equal.
\end{defi}

\begin{defi}\label{def:tassel}
We say that the family $T$ is a {\em tassel} based at $b$ if
\begin{enumerate}
\item $\range(T)$ is contained in a contractible open subset of $M$.
\item $\alpha_i(0)=b$ for $i=1,\ldots,n$.
\item there is a parameterization $(x_1,\ldots,x_d)\mapsto M$ such that
      $b=(0,\ldots,0)$ and $\alpha_i$ can be written as a graph
      $\alpha_i(t)=(t,f_i(t))$, where $f_i:[0,t_i] \to \R^{d-1}$.
\item if there is a regular point in $\range(T)$ with a certain type, then
      there are points with the same type in every neighborhood of $b$.
\item all the curves in the family are different.
\end{enumerate}
\end{defi}

\begin{defi}
A family of curves
$\alpha_{1,1},\ldots,\alpha_{1,n_1},\ldots,\alpha_{k,1},\ldots,\alpha_{k,n_k}$
is a {\em web} if the curves $\alpha_{j,1}\dots\alpha_{j,n_j}$ form a tassel
$T_j$ for $j=1,\ldots,k$ and curves in two different tassels do not intersect
(except, possibly, at their endpoints).
\end{defi}

In~\cite{BS1} is proven the following proposition.

\begin{prop}
For every family of curves $F$ in $M$ there is a web $W$ such that all curves
in $F$ are products of curves in $W$ (or their inverses).
\end{prop}

Applying this proposition to a loop $\gamma$ we obtain a web $W$ formed by a
family of curves $\alpha_1,\ldots,\alpha_n$ such that $\gamma$ can be obtained
as product of curves in $W$, then repeating the same construction as in proof
of theorem~\ref{teo:decomp-an-loop} we can construct loops
$\gamma_i=\beta(e_-(\alpha_i))\alpha_i\beta(e_+(\alpha_i))$ then $\gamma$ is
a product of loops $\gamma_i$.

\begin{defi}
Let $T$ be a tassel, $G_T$ is the closed subgroup of $G^n$ generated by all
the types of regular points in $T$.
\end{defi}

Now let us see what are the possible holonomies for curves in a tassel,
let $T$ be a tassel composed by curves $\alpha_1,\ldots,\alpha_n$, $\range(T)$
is contained in a contractible open set $U$, then we can fix a trivialization
of the bundle over $U$ and associate to each connection a element of $G$ over
each curve, then we can identify the set of possible values of holonomies for
the curves $\alpha_i$ with a subset of $G^n$, in~\cite{BS1}
is also proved:

\begin{prop}
For a tassel $T$ the set of possible values for the holonomies is $G_T$
\end{prop}

\begin{lem}
If $G$ is semisimple then $G_T=G^n$.
\end{lem}

\Proof
It is sufficient to show that
all the elements of $G^n$ of the form $E(g,i)=(e,\ldots,g,\ldots,e)$ (the
element of $G^n$ that has $g$ in the $i$-th component and $e$ in the
others) are in $G_T$, now we say that a element of $G^n$ is of the form
$E(g,i,i_1,\ldots,i_k)$ if its $i$-th component is $g$ and the components
$i_1,\ldots,i_k$ are $e$. We will show that $G_T$ has elements of the form
$E(g,i,i_1,\ldots,i_k)$ for every $g$, $i$, $i_1,\ldots,i_k$ (when $k=n-1$
this implies that $G_T$ has all the elements of the form $E(g,i)$ then
$G_T=G^n$), we proceed by induction in $k$, when $k=1$, given $g$, $i$,
$i_1$, we can find a regular point of $T$ such that $\alpha_i$ passes
trough $p$ and $\alpha_{i_1}$ does not, then the type of $p$ is a element of
$G_T$ of the form $E(g,i,i_1)$.
To proceed with the inductive step assume that we are given
$g=g_1g_2g_1^{-1}g_2^{-1}$, $i$, $i_1,\ldots,i_{k+1}$.
Then take elements of $G_T$ $\hat{g}_1$ of the form $E(g_1,i,i_1,\ldots,i_k)$
and $\hat{g}_2$ of the form $E(g_2,i,i_2,\ldots,i_{k+1})$ (they exist by
induction hypothesis), then $\hat{g}_1\hat{g}_2\hat{g}_1^{-1}\hat{g}_2^{-1}$
is a element of $G_T$ of the form $E(g,i,i_1,\ldots,i_{k+1})$.
Since $G$ is semisimple elements of the form $g_1g_2g_1^{-1}g_2^{-1}$ generate
$G$, then for every $g \in G$ $G_T$ has elements of the form 
$E(g,i,i_1,\ldots,i_{k+1})$.
\QED

Now we can give a proof of theorem~\ref{teo:hg-diff-loop-nonsolv}

\ProofP{of~\ref{teo:hg-diff-loop-nonsolv}}
Let $(E,G,M,\Pi)$ be a principal bundle with $G$ a non solvable group, let
$\hat{G}$ be the quotient of $G$ by its radical, then $\hat{G}$ is
semisimple, let $(\hat{E},\hat{G},M,\hat{\Pi})$ be the extension of the
bundle to the group $\hat{G}$.
Let $\gamma$ be a loop in $M$ not null in $\loopg$ then as remarked before we
can decompose $\gamma$ in product of loops
$\gamma_i=\beta(e_-(\alpha_i))\alpha_i\beta(e_+(\alpha_i)$ where the curves
$\alpha_1,\ldots,\alpha_n$ form a web, then because $\hat{G}$ is semisimple
we can choose holonomies independently for the loops $\gamma_i$ then we can
proceed as in the proof of theorem~\ref{teo:hg-an-loop} and find a connection
$A$ in $\hat{E}$ such that $H_A(\gamma) \ne e$ then we can pull-back this
connection to $E$ and obtain a connection in $E$ such that the holonomy of
$\gamma$ is not null.
\QED

Now let us prove the theorem~\ref{teo:hg-diff-loop}.
It follows from the following proposition.

\begin{prop}
Let $\gamma$ be a $C^N$ loop such that $H_A(\gamma)=e$ for all connections in
the bundle $(E,G,M,\Pi)$ then there are loops in $E_G(\loopg)$ arbitrarily
$C^N$ close to $\gamma$.
\end{prop}

\Proof
We first decompose $\gamma$ in curves forming a web with tassels 
$T_1,\ldots,T_n$, we will do small deformations to these curves to obtain a
family of curves that intersects only in a finite number of segments or
isolated points.
We need to do such deformation in a way that the holonomy of the deformed loop
$\hat{\gamma}$ be $e$ for every connection, then the same argument as in proof
of theorem~\ref{teo:hg-an-loop} shows that $\hat{\gamma} \in E_G(\loopg)$.
To accomplish this is sufficient that the group of possible holonomies of the
deformed tassel $G_{\hat{T_i}}$ be included in $G_{T_i}$, and for this is 
sufficient that the deformation do not take apart curves that intersect.

Fix $\epsilon>0$, take $\alpha_1,\dots,\alpha_c$ curves in a tassel $T$ in the
decomposition of $\gamma$ and take the parameterization in
definition~\ref{def:tassel} such that $\alpha_1(t)=(t,0,\ldots,0)$,
$t\in [0,1]$, we identify points in $M$ with the corresponding points in
$\R^d$ in the parameterizations and identify points $t$ in $[0,1]$ with
points $(t,0,\ldots,0)$. Also take a $C^\infty$ function $\rho:\R \to [0,1]$
such that $\rho(x)=1$ if $x \le -1$ or $x \ge 1$ and $\rho(x)=0$ if
$-1/2 \le x \le 1/2$, let $a_n=\max_{x\in\R}|\rho^{(n)}(x)|$ and let
$\rho_{p,\delta}(x)=\rho((x-p)/\delta)$, note that
$|\rho_{p,\delta}^{(n)}(x)| \le a_n/\delta^n$.

We say that a point in the intersection of two curves is singular if it is
not in the interior of a common interval of both curves. Let $A$ be the set
of singular points of intersection between $\alpha_1$ and the other curves,
because of the way in what was chosen the coordinate system points in $A$
are of the form $(t,0,\ldots,0)$ and thus we identify them with points of
$[0,1]$ which are the values of the parameter of the point taking the
parameterization $\alpha_i(t)=(t,f_i(t))$ (where $f_i$ is as in
definition~\ref{def:tassel}).
Let $A'$ be the set of accumulation points of $A$ thus $A'$ is a compact
subset of $[0,1]$.

Let $p \in A'$, we define $\cC^N(p)$ as the set of curves $\alpha_i$ that
have a contact of order $N$ with $\alpha_1$ in $p$ (that is the first $N$
derivatives of $f_i$ in $p$ are null), note that if $p$ is accumulation
point of intersection points of $\alpha_1$ and $\alpha_i$ then 
$\alpha_i \in \cC^N(p)$.
And define $\cC^{<N}(p)$ the set of
curves $\alpha_i$ having a contact of order less than $N$ with $\alpha_1$ in
$p$. We take a $\delta=\delta_p$ such that
\begin{enumerate}
\item $\delta<1$ and if $p \ne t_i$ then $t_i \not\in (p-\delta,p+\delta)$
      ($t_i$ as in definition~\ref{def:tassel}).
\item the curves in $\cC^N(p)$ don't intersect the curves in $\cC^{<N}$ in
      other point than $p$ for parameters values in $(p-\delta, p+\delta)$.
\item for each $\alpha_i \in \cC^N(p)$, $f_i$ is a $C^N$ function having
      the first $N$ derivates equal to zero, then there exist
      $r_p(f_i^{(n)},\delta)$ such that
\[
        |f_i^{(n)}(x)|<r_p(f_i^{(n)},\delta)|x-p|^{N-n}
\]
      for all $n \le N$, then take
      $\delta$ such that $r_p(f_i^{(n)},\delta)<\epsilon$ and
\[
        a_n r_p(f_i,\delta) + \cdots + C^n_k a_k r_p(f_i^{(n-k)},\delta) +
            \cdots + r_p(f_i^{(n)},\delta) < 2r_p(f_i^{(n)},\delta).
\]
      where $C^n_k$ is the binomial coefficient, note that this implies that
      $|f_i^{(n)}(x)|<\epsilon$ for $x\in(p-\delta,p+\delta)$, $n \le N$.
\end{enumerate}

Take $P=\{p_1,\ldots,p_m\} \subset A'$ such that
$(p_i-\delta_i/2,p_i+\delta_i/2)$ is a finite cover of $A'$.

For each $\alpha_i$ take $Q_i \subset P$
such that $q \in Q_i$ if and only if $\alpha_i \in \cC^N(q)$,
for $x \in \R$ let
\begin{eqnalign*}
    Q_i^-(x)&=\{q \in Q_i: q<x \hbox{ and there is no $q<q'<x$ such that }
                \cC^N(q)=\cC^N(q')\}  \\
    Q_i^+(x)&=\{q \in Q_i: x<q \hbox{ and there is no $x<q'<q$ such that }
                \cC^N(q)=\cC^N(q')\}
\end{eqnalign*}
and let $Q_i(x)=Q_i^-(x) \cup Q_i^+(x)$ note that $\# Q_i^\pm(x) < 2^c$,
next define
\[
    \rho_i(x)= \prod_{q \in Q_i(x)} \rho_{q,\delta_q}(x)
\]
let $\bar{f}_i=f_i\rho_i$ and $\bar{\alpha}_i(t)=(t,\bar{f}_i(t))$.
Let $x \in [0,t_i]$, let
\[
    \{q_1,\ldots,q_k\}=\{q \in Q_i(x) : x \in (q-\delta_q, q+\delta_q) \}
\]
and $\delta_j=\delta_{q_j}$ ordered such that
$\delta_1>\delta_2>\cdots>\delta_k$, let ${}_0f_i=f_i$,
${}_jf_i={}_{j-1}f_i \rho_{q_j,\delta_j}$ then $\bar{f}_i(x)={}_kf_i(x)$ we
will prove by induction
$|{}_jf_i^{(n)}(x)| < 2^j r_{q_k}(f^{(n)},\delta_j)|x-q_k|^{N-n}$, for $j=0$
is clear, for $j>0$
\begin{eqnalign*}
    |{}_jf_i^{(n)}(x)| &= \left| \sum_{\ell=0}^n C^n_\ell
                   \rho_{q_j,\delta_j}^{(\ell)}(x) {}_{j-1}f_i^{(n-\ell)}(x)
                   \right|  \\
        &\le \sum_{\ell=0}^n C^n_\ell \left| \frac{a_\ell}{\delta_j^\ell}
        2^{j-1} r_{q_k}(f_i^{(n-\ell)},\delta_j)|x-q_k|^{N-n+\ell} \right| \\
        &< 2^j r_{q_k}(f_i^{(n)},\delta_j)|x-q_k|^{N-n}
\end{eqnalign*}
where its used that $|x-q_k|<\delta_k<\delta_j$.
Then 
\[
   |\bar{f}^{(n)}_i(x)| < 2^k r(f_i^{(n)},\delta_k)|x-q_k|^{N-n} < 
                         2^k\epsilon<2^{2^{c+1}}\epsilon.
\]

This shows that the distance of $\alpha_i$ and $\bar{\alpha}_i$ is less than
$2^{2^{c+1}}\epsilon$ in the $C^N$ topology, and this construction removed
all the accumulation points of singular intersection points between $\alpha_1$
and $\bar{\alpha}_2$, then $\alpha_1$ and $\bar{\alpha}_2$ intersect in a
finite number of intervals or isolated points.
We have to see that $G_{\bar{T}}$, the set of possible values for holonomies
of the curves $\bar{\alpha}_i$, is not larger than $G_T$, to show this it is
sufficient to show that if $\alpha_i(t)=\alpha_j(t)$ then
$\bar{\alpha}_i(t)=\bar{\alpha}_j(t)$.

Let $Q'_i(x)=\{q \in Q_i(x) : x \in (q-\delta_q,q+\delta_q)\}$, then
\[
    \rho_i(x)= \prod_{q \in Q'_i(x)} \rho_{q,\delta_q}(x)
\]
if $q \in Q_i(x)$ and $q \not\in Q_j(x)$ then $\alpha_i \in \cC^N(q)$ and
$\alpha_j \in \cC^{<N}(q)$ thus $\alpha_i$ do not intersect $\alpha_j$ in
$(q-\delta_q,q+\delta_q)$ and then $\alpha_i(x) \ne \alpha_j(x)$.
This shows that if $\alpha_i(x)=\alpha_j(x)$ then $Q'_i(x)=Q'_j(x)$ then
$\rho_i(x)=\rho_j(x)$ and this implies
$\bar{\alpha}_i(x)=\bar{\alpha}_j(x)$.
\QED

\end{document}